\documentclass[a4, 12pt]{article}

\usepackage{mathptmx}
\usepackage{amsmath}
\usepackage{amssymb}
\usepackage{enumerate}
\usepackage[margin=2cm]{geometry}

\newtheorem{theorem}{Theorem}

\newtheorem{lemma}{Lemma}
\newtheorem{definition}{Definition}
\newtheorem{remark}{Remark}
\newcommand{\qed}{$\mbox{}\hfill\Box$}

\begin{document}

\title{New bounds of extended energy of a graph}

\author{Abujafar Mandal\thanks{Email: abujafarmondal786@gmail.com} \and Sk. Md. Abu Nayeem\thanks{Corresponding author. Email: nayeem.math@aliah.ac.in}}
\date{}






\maketitle

\begin{abstract} 
The extended adjacency matrix of a graph with $n$ vertices is a real symmetric matrix of order $n\times n$ whose $(i,j)$-th entry is the average of the ratio of the degree of the vertex $i$ to that of the vertex $j$ and its reciprocal when $i,j$ are adjacent and zero otherwise. The aggregate of absolute eigenvalues of the extended adjacency matrix is termed the extended energy. 

In this paper, the concept of extended vertex energy is introduced, and some bounds of extended vertex energy are obtained. From there, we establish some new upper bounds of the extended energy of a graph involving order, size, largest, and smallest degree. We show that those are improvements of some existing bounds. Through direct manipulation, we have also established some more upper and lower bounds of extended energy, which are either better or incomparable with the existing bounds. Finally, some improved bounds of Nordhaus-Gaddum-type are found.

\medskip

\noindent MSC (2020): 05C50, 05C07.

\medskip
\noindent
\textit{Keywords.} Extended energy; extended vertex energy; Nordhaus-Gaddum-type bound.
\end{abstract}

\section{Introduction}
For a graph $G$, let $V(G)$ represent its vertex set and $E(G)$ its edge set. If the order and size of $G$ are $|V(G)|=n$ and $|E(G)|=e$ respectively, then $G$ is said to be an $(n,e)$-graph. By a simple graph, we mean a graph with no parallel edges and no self-loops. In this work, our concern is limited to simple graphs only. Let $d_i$ denote the degree of $i\in V(G)$; the smallest and largest among all $d_i$'s be $\delta_m$ and $\delta_M$ respectively. If $d_i=r$ for all $i\in V(G)$, then $G$ is an $r$-regular graph. Evidently, the complete graph $K_n$ is $\hat{n}$-regular, where $\hat{n}=n-1$. For convenience, we use $i\sim j$ to denote that $i$ and $j$ are adjacent vertices. 

For a simple graph $G$ of order $n$, the $n\times n$ symmetric matrix whose $ (i,j)$-th entry is 1 if $i\sim j$ and 0 otherwise is called the adjacency matrix of $G$, and it is denoted by $A(G)$. Since it is real and symmetric, all eigenvalues are real. The eigenvalues $\lambda_1,\lambda_2,\ldots,\lambda_n$ of $A(G)$, indexed in the decreasing order of their values, are also known as the eigenvalues of $G$. Gutman and co-researchers \cite{gr75,gt72} found that the $\pi$-electron energy of a molecule can be estimated by the aggregate of absolute eigenvalues of the corresponding molecular graph $G$. Hence that sum was termed graph energy (see \cite{g01,glz09,gsl12}). Thus, $\varepsilon(G)$, denoting the energy of $G$, is defined by $\varepsilon(G)=\sum_{j=1}^{n}|\lambda_j|$. Graph energy and several variations of it have been studied by many researchers \cite{f22,g01,glz09,gr75,gsl12,gt72,km01} in the recent past.
	
In 1994, Yang et al. \cite{yx94} brought up a generalization of the adjacency matrix which they called extended adjacency matrix. For a graph $G$ of order $n$, the extended adjacency matrix is an $n\times n$ symmetric real matrix and it is denoted by $A_{ex}(G)$. The $(i,j)$-th entry of $A_{ex}(G)$ is $\frac{1}{2}(\frac{d_i}{d_j}+\frac{d_j}{d_i})$ if $i\sim j$, and $0$, if $i,j$ are non-adjacent. Thus, $A_{ex}(G)$ turns down to $A(G)$ when $G$ is regular. The extended adjacency index \cite{yx94} of a graph $G$ is defined as $$EA(G)=\frac{1}{2}\sum_{i\sim j}\left(\frac{d_i}{d_j}+\frac{d_j}{d_i}\right).$$
All the eigenvalues $A_{ex}(G)$ are also real. Let the eigenvalues of $A_{ex}(G)$ be $\eta_1,\eta_2,\ldots,\eta_n$ indexed in the decreasing order of their values.  Yang et al. \cite{yx94} put forward the notion of the extended energy of a graph $G$ as $$\varepsilon_{ex}(G)=\sum_{j=1}^{n}\lvert \eta_j\rvert.$$ 

Adiga and Rakshith \cite{ar18}, Das et al. \cite{dg17}, Liu et al. \cite{lp20}, Ghorbani et al. \cite{ga22,gl21}, Wang et al. \cite{wm19} have obtained several bounds of extended energy of a graph. Gutman \cite{g17} has shown that for a bipartite graph, the extended energy is not less than the ordinary energy, and has conjectured that it holds for any graph, in general.

In 2018, the notion of vertex energy of a graph was put forward by Arizmendi et al. \cite{oj18}. In the present work, we introduce the idea of extended vertex energy and find some lower and upper bounds of it. This allows us to derive two new upper bounds for the extended energy of a graph. We show that those bounds are improvements of some existing bounds reported in \cite{dg17, ga22} and \cite{wm19}. Through direct manipulation, we have also established some more upper and lower bounds of extended energy, which are either better or incomparable with the existing bounds given in \cite{ar18, dg17, ga22, lp20} and \cite{wm19}.
	
The last section of the paper contains some bounds of Nordhaus-Gaddum-type for extended energy and it is shown that those are also improvements of similar bounds obtained in \cite{wm19}.

\section{Extended energy of a vertex and its bounds}
	
We denote the trace of a square matrix $M$ by $Tr(M)$. Also, we denote the absolute value of a real square matrix $M$ by $|M|$, i.e., $|M|=(MM^T)^{\frac{1}{2}}.$ For ease of notation, we denote the adjacency matrix of a graph $G$ by $A$ and the extended adjacency matrix of $G$ by $A_{ex}$ simply. 
Nikiforov \cite{nik07} showed that the Schatten 1-norm of $A$ gives the energy of $G$, i.e., $\varepsilon(G)=Tr(|A|)$. Keeping this in mind, in a recent development towards the theory of graph energy, the notion of vertex energy of a graph has been introduced in \cite{oj18}, by Arizmendi et al. The energy of the $i$-th vertex, denoted by $\varepsilon_{i}(G)$ or simply by $\varepsilon_i$, is given by 
	$$\varepsilon_i=|A|_{ii}\mbox{ for } i\in V(G)$$ and as such $$\varepsilon(G)=\varepsilon_1+\varepsilon_2+\cdots+\varepsilon_n.$$

In a similar fashion, we introduce here the idea of extended vertex energy of a graph. We denote the extended energy of $i\in V(G)$ by $\varepsilon_{ex_i}(G)$ or simply by $\varepsilon_{ex_i}$ and define as $$\varepsilon_{ex_i}=|A_{ex}|_{ii} \mbox{ for } i\in V(G)$$  and so $$\varepsilon_{ex}(G)=\varepsilon_{ex_1}+\varepsilon_{ex_2}+\cdots+\varepsilon_{ex_n}.$$

In fact, the extended energy of a vertex quantifies the interaction of that vertex with other vertices of the graph as its contribution towards the extended energy of the whole graph. It is easy to realize that the extended vertex energy is not impacted by the vertices which belong to some other component of the graph, or more precisely $\varepsilon_{ex_i}(G)=\varepsilon_{ex_i}(G_i)$, where the $i$-th vertex is contained in the component $G_i$ of $G$.

Concisely, we give the definition of the extended vertex energy as follows.
	
	\begin{definition} Let the vertices of a graph $G$ be $1,2,\ldots,n$ and $A_{ex}$ denotes the extended adjacency matrix of $G$. For each $i\in V(G)$, the extended energy of the $i$-th vertex of $G$, denoted by $\varepsilon_{ex_{i}}(G)$ or more conveniently by $\varepsilon_{ex_i}$, is given by $$\varepsilon_{ex_i}=|A_{ex}|_{ii}$$
	where $|A_{ex}|=\left(A_{ex}A_{ex}^T\right)^\frac{1}{2}=\left(A_{ex}^2\right)^\frac{1}{2}$.
	\end{definition}

The following lemma helps us to obtain some bounds of extended vertex energy.

\begin{lemma}\label{exvarenergyre}
		Let $G=(V,E)$ be an $n$ vertices graph where the vertices are labelled as $1,2,\ldots,n$. Then for each $i\in V(G)$
		\begin{eqnarray}
			\varepsilon_{ex_i}=\sum_{r=1}^{n}q_{ir}|\eta_r|
		\end{eqnarray}
	where $\eta_r$ is the $r$-th eigenvalue of $A_{ex}$, the extended adjacency matrix of the graph $G$ and the weights $q_{ir}$ satisfy $$\sum_{i=1}^{n}q_{ir}=1 \mbox{ for each } r\in V(G) \mbox{ and } \sum_{r=1}^{n}q_{ir}=1 \mbox{ for each } i\in V(G).$$  
	Moreover, $q_{ir}=w_{ir}^2$ where $W=(w_{ir})$ is an orthogonal matrix formed by the eigenvectors of $A_{ex}$ as column vectors.
	\end{lemma}

\noindent\textit{Proof.} Since $A_{ex}$ is a real symmetric matrix, it is orthogonally diagonalizable. Let $\{w_1, w_2,\ldots, w_n\}$ be the orthonormal set of eigenvectors of $A_{ex}$, with corresponding eigenvalues $\eta_1, \eta_2,\ldots, \eta_n$. If $W=(w_{ir})$ be the matrix made by eigenvectors, i.e., the columns of $W$ are $w_1, w_2,\ldots, w_n$, and $H$ be the diagonal matrix $\operatorname{diag}(\eta_1, \eta_2,\ldots, \eta_n)$, then $W$ is an orthogonal matrix and $A_{ex}=WHW^T.$

Thus
\begin{eqnarray}
	|A_{ex}|=\left(A_{ex}^2\right)^{\frac{1}{2}}=W|H|W^T,
\end{eqnarray}
where $|H|=\operatorname{diag}(|\eta_1|, |\eta_2|,\ldots, |\eta_n|)$.

Since $|H|_{ir}=|\eta_i|$ if $i=r$, and 0, otherwise, 
\begin{equation}
	|A_{ex}|_{ii}=\sum_{r,h=1}^{n}W_{ih}|H|_{hr}W_{ri}^T=\sum_{r=1}^{n}w_{ir}|\eta_r|w_{ir} =\sum_{r=1}^{n}w_{ir}^2|\eta_r|=\sum_{r=1}^{n}q_{ir}|\eta_r|,
\end{equation}
which is our desired result.
\qed

\begin{lemma}[McClelland] {\cite{m71}} \label{mcbd}
	Let $G$ be an $(n,e)$-graph. Then $\varepsilon(G) \leq \sqrt{2ne}$
	with equality if and only if $G \cong \frac{n}{2}K_2$ when $n$ is even, or $G\cong \overline{K_n}$ .
\end{lemma}

\begin{lemma}\cite{ga22}
    Let $G$ be an $n $  vertices graph and $EA(G)$ is the extended adjacency index. If $\varepsilon_{ex}(G)$ is the extended energy of $G$, then $$\varepsilon_{ex}(G)\leq \sqrt{2n}EA(G).$$
\end{lemma}

\subsection{Upper bounds of extended vertex energy}
In this subsection, we find two upper bounds of extended vertex energy. 

\begin{theorem} \label{ubexvertex}
	Let the vertices of a graph $G$ be $1,2,\ldots,n$. If the smallest degree and the largest degree of $G$ are $\delta_m$ and $\delta_M$ respectively and $d_i$ be the degree of the $i$-th vertex, then $$\varepsilon_{ex_i}\leq \frac{1}{2}\left(\frac{\delta_M}{\delta_m}+\frac{\delta_m}{\delta_M}\right)\sqrt{d_i}$$
	and the sufficient and necessary condition for equality is that the component consisting of the vertex $i$ is isomorphic to a star graph with $i$ being its center.
\end{theorem}

\noindent\textit{Proof.} 
	Let us define a positive linear functional $\psi_i:M_n(\mathbb{R}) \to \mathbb{R}$ by $$\psi_i(X)\mapsto (X)_{ii}.$$
	
	Then $$Tr(X)=\psi_1(X)+\psi_2(X)+\cdots+\psi_n(X).$$
	
	Thus 
	\begin{equation}
	\varepsilon_{ex}(G)=Tr(|A_{ex}|)=\psi_1(|A_{ex}|)+\psi_2(|A_{ex}|+\cdots+\psi_n(|A_{ex}|).	
\end{equation}

Since all $\psi_i,i\in V(G)$ are positive linear functionals on $M_n(\mathbb{R})$, they satisfy the H\"{o}lder's inequality. Thus for $X,Y\in M_n(\mathbb{R})$ and positive real numbers $p,q$ satisfying $\frac{1}{p}+\frac{1}{q}=1$, we have
	\begin{eqnarray}
		\psi_i(|XY|)\leq \psi_i(|X|^p)^{\frac{1}{p}}\psi_i(|Y|^q)^{\frac{1}{q}}.
	\end{eqnarray}
	
Choosing $X=A_{ex}, Y=I_n,p=2, q=2$ and then squaring both sides, we get
\begin{equation}\psi_i(|A_{ex}|)^2 \leq \psi_i(|A_{ex}|^2)=\psi_i(A_{ex}^2).\label{ineq1}\end{equation}

Now, 
\begin{eqnarray}{\nonumber}
	\psi_i(A_{ex}^2)&=&\sum_{i\sim r}\frac{1}{4}\left(\frac{d_i}{d_r}+\frac{d_r}{d_i}\right)^2 \nonumber	\\
	&\le&\frac{1}{4}\sum_{i\sim r}\left(\frac{\delta_M}{\delta_m}+\frac{\delta_m}{\delta_M}\right)^2\label{ineq2} =\frac{1}{4}\left(\frac{\delta_M}{\delta_m}+\frac{\delta_m}{\delta_M}\right)^2d_i. 
\end{eqnarray}

Hence, $$\psi_i(|A_{ex}|)^2 \leq \frac{1}{4}\left(\frac{\delta_M}{\delta_m}+\frac{\delta_m}{\delta_M}\right)^2d_i.$$

Taking positive square roots on both sides, we get $$\varepsilon_{ex_i} \leq \frac{1}{2}\left(\frac{\delta_M}{\delta_m}+\frac{\delta_m}{\delta_M}\right)\sqrt{d_i}.$$

It is uncomplicated to observe that the above relation is an equality if $i$ is the centre of a star. 
For equality in (\ref{ineq2}), it is required that $d_i=\delta_M$ and $d_j=\delta_m$ for all $j\sim i$. For equality in (\ref{ineq1}), it is required that the distribution of the matrix $A_{ex}$ to be Dirac measure $\delta_x$ for some $x$ (see \cite{oj18}). Thus the sufficient and necessary condition for equality being the component, to which the vertex $i$ belongs, is isomorphic to a star graph with $i$ being its center.

\begin{theorem}\label{ubdvertexfth}
	In a graph $G$ with $n$ vertices, $e$ edges and smallest degree $\delta_m$, for each $i\in V(G)$ $$\varepsilon_{ex_i} \leq \sqrt{ \frac{1}{4\delta_m^2}\sum_{i\sim j}d_j^2+\frac{d_i^3}{4\delta_m^2}+\frac{1}{2}d_i}.$$
\end{theorem}

\noindent\textit{Proof.} As before, let us consider a linear positive functional $\psi_i:M_n(\mathbb{R}) \to \mathbb{R}$ as $$\psi_i(X)\mapsto (X)_{ii}.$$

Then,
$\psi_i(|A_{ex}|)^2 \leq \psi_i(|A_{ex}|^2)=\psi_i(A_{ex}^2).$

Now, 
\begin{eqnarray}
	\psi_i(A_{ex}^2)	
	&=&\sum_{i\sim j}\frac{1}{4}\left(\frac{d_i}{d_j}+\frac{d_j}{d_i}\right)^2\nonumber\\
	&=&\frac{1}{4}\sum_{i\sim j}\left(\frac{d_i^2}{d_j^2}+\frac{d_j^2}{d_i^2}\right)+\frac{1}{2}d_i \nonumber\\
	&=&\frac{d_i^2}{4}\sum_{i\sim j}\frac{1}{d_j^2}+\frac{1}{4d_i^2}\sum_{i\sim j}d_j^2+\frac{1}{2}d_i\nonumber \\
	&\leq&\frac{d_i^2}{4}\sum_{i\sim j}\frac{1}{\delta_m^2}+\frac{1}{4\delta_m^2}\sum_{i\sim j}d_j^2+\frac{1}{2}d_i\nonumber \\
	&\leq&\frac{d_i^3}{4\delta_m^2}+\frac{1}{4\delta_m^2}\sum_{i\sim j}d_j^2+\frac{1}{2}d_i.\nonumber
\end{eqnarray}

Since $\psi_i(|A_{ex}|)=\varepsilon_{ex_i}$, we have $$\varepsilon_{ex_i}^2 \leq \frac{1}{4\delta_m^2}\sum_{i\sim j}d_j^2+\frac{d_i^3}{4\delta_m^2}+\frac{1}{2}d_i.$$ 

Taking positive square roots on both sides, we get $$\varepsilon_{ex_i} \leq \sqrt{ \frac{1}{4\delta_m^2}\sum_{i\sim j}d_j^2+\frac{d_i^3}{4\delta_m^2}+\frac{1}{2}d_i}.$$
\qed

\bigskip

\subsection{A lower bound of extended vertex energy}
Now we obtain a simple lower bound of extended vertex energy $\varepsilon_{ex}(v_i)$. We shall use the following lemma for that purpose.

\begin{lemma}\cite{gl21}\label{lbdsprth}
	Let $G$ be a graph with $\lambda_1$ and $ \eta_1$ being spectral radii of $A$ and $A_{ex}$ respectively. Then, $\lambda_1\leq \eta_1\leq \frac{1}{2}\left(\frac{\delta_M}{\delta_m}+\frac{\delta_m}{\delta_M}\right)\lambda_1$. The sufficient and necessary condition for equality of the left inequality is that $G$ is a regular graph and that for the right inequality is that $G$ is a complete bipartite graph.
\end{lemma}

\begin{theorem}\label{lbdvexe01}
	Let $G$ be an $(n,e)$-graph, where $e>0$. Then for each $i\in V(G)$ $$\varepsilon_{ex_i}\geq \frac{d_i}{k}$$ where $k=\frac{1}{2}\left(\frac{\delta_M}{\delta_m}+\frac{\delta_m}{\delta_M}\right)\delta_M$ and the sufficient and necessary condition for equality being $G \cong K_{\delta_M,\delta_M}$.
	\end{theorem}

\noindent\textit{Proof.} 
From Lemma \ref{lbdsprth}, we have $\eta_1\leq \frac{1}{2}\left(\frac{\delta_M}{\delta_m}+\frac{\delta_m}{\delta_M}\right)\lambda_1\le \frac{1}{2}\left(\frac{\delta_M}{\delta_m}+\frac{\delta_m}{\delta_M}\right)\delta_M=k$. It can be easily verified that $|x|\geq x^2$ for  $x\in[-1,1]$. So, $$\left|\frac{x}{k}\right|\geq \left(\frac{x}{k}\right)^2 \mbox{ for } x\in[-k,k]$$ with equality holds only if $x \in \{-k,0,k\}$.

From Lemma \ref{exvarenergyre}, we have $$\varepsilon_{ex_i}=\sum_{j=1}^{n}q_{ij}|\eta_j|\geq \sum_{j=1}^{n}q_{ij}\frac{\eta_j^2}{k} \geq \frac{d_i}{k}\cdot$$

Now $k$ and $-k$ will be eigenvalues of $A$ if $G$ is regular and  bipartite (see  \cite{oj18}). Again, if $G$ is regular, then $A_{ex}=A$ and so $\eta_i=\lambda_i$ and $k=\delta_M$.  Thus $\eta_i \in \{-k, 0, k\}$ if and only if $G$ is isomorphic to a regular bipartite graph.
 \qed

\section{Upper bounds of extended energy}

For a graph $G$, the forgotten index, denoted by $F(G)$ or simply by $F$, is defined by $$F(G)=\sum_{r=1}^{n}d_r^3=\sum_{r\sim i}\left(d_r^2+d_i^2\right).$$ The modified second Zagreb index is defined as $$M_2^*(G)=\sum_{r\sim i}\frac{1}{d_rd_i}\cdot$$

The following upper bounds of extended energy were obtained by Adiga et al.
\begin{lemma}\cite{ar18} \label{ubdexadiga}
    Let $G$ be a graph of order $n$. Then $$\varepsilon_{ex}(G)\leq \frac{\delta_M}{\delta_m}\varepsilon(G).$$
\end{lemma}

The following upper bounds of extended energy were obtained by Das et al.

\begin{lemma}\cite{dg17}\label{dg17ubd}
	Let $G$ be a $(n,e)$-graph. If the largest degree is $\delta_M$ and the smallest degree is $\delta_m$, then 
\begin{enumerate}\item[(i)] $\displaystyle \varepsilon_{ex}(G)\leq \frac{1}{2}\left(\frac{\delta_M}{\delta_m}+\frac{\delta_m}{\delta_M}\right)\sqrt{2ne}$, 
\item[(ii)] $\displaystyle \varepsilon_{ex}(G)\leq \sqrt{\left(\frac{\delta_M}{\delta_m}+\frac{\delta_m}{\delta_M}\right)}\sqrt{\frac{nF}{2\delta_m^2}}$
and the both equality holds if and only if $G\cong \frac{n}{2}K_2$ for even $n$.
\end{enumerate}
\end{lemma}

The following upper bound is given by Liu et al. \cite{lp20}.
\begin{lemma}\cite{lp20}\label{ubexlm05}
    Let G be a non-empty graph of order n with m edges. Then
\begin{enumerate}[(a)]
    \item $$\varepsilon_{ex}(G)\leq \tau +\sqrt{(n-1)(2M-2\tau^2)}$$
    where $\tau =\max\left\{\sqrt{\frac{2M}{n}},\frac{4\delta_m^2M_2^*-2e}{n}, \frac{EA(G)}{n}, \frac{2e}{n}\right\}$ and $2M=\sum_{i=1}^{n}\eta_i^2$. Equality holds if and only if $G\cong K_n$ or $G\cong \frac{n}{2}K_2$.
    \item $$\varepsilon_{ex}(G)\leq \tau +\sqrt{(n-1)\left(\frac{1}{2}\left(\frac{\delta_M}{\delta_m}+\frac{\delta_m}{\delta_M}\right)EA(G)-2\tau^2\right)}$$
    where $\tau =\max\left\{\sqrt{\frac{1}{2}\left(\frac{\delta_M}{\delta_m}+\frac{\delta_m}{\delta_M}\right)\frac{EA(G)}{n}},\frac{4\delta_m^2M_2^*-2e}{n}, \frac{EA(G)}{n}, \frac{2e}{n}\right\}$. Equality holds if and only if $G\cong K_{\frac{n}{2},\frac{n}{2}}$ or $G\cong \frac{n}{2}K_2$.
    
\end{enumerate}

Using the upper bounds of extended vertex energy given in Theorem \ref{ubexvertex} and Theorem \ref{ubdvertexfth}, we obtain two new upper bounds of $\varepsilon_{ex}(G)$. For that purpose, the lemma that follows will be applied.

\end{lemma}

\begin{lemma}\label{poinq}
	Let $t_1, t_2,\ldots, t_r$ be $r$ nonnegative real numbers. Then $$\sum_{i=1}^r\sqrt{t_i}\leq \sqrt{r\sum_{i=1}^rt_i}.$$
	The necessary and sufficient condition for equality is that $t_1=t_2=\cdots=t_r$.
\end{lemma}

\noindent\textit{Proof.} 
	Since  $t_1, t_2,\ldots, t_r$ are nonnegative real numbers, so,  $\sqrt{t_1}, \sqrt{t_2},\ldots, \sqrt{t_r}$ are also nonnegative real numbers. Now using AM-QM inequality \cite{bm88} on $\sqrt{t_1}, \sqrt{t_2},\ldots, \sqrt{t_r}$ we have the required result.
\qed

\begin{theorem}\label{ubtexenergy}
	Let $G$ be an $(n,e)$-graph. If the largest degree and the smallest degree of $G$ are $\delta_M$ and $\delta_m$ respectively, then
	$$\varepsilon_{ex}(G)\leq \frac{1}{2}\left(\frac{\delta_M}{\delta_m}+\frac{\delta_m}{\delta_M}\right)\left(\sqrt{(n-2)(2e-\delta_m-\delta_M)}+\sqrt{\delta_m}+\sqrt{\delta_M}\right).$$ The sufficient and necessary condition for the relation to hold with equality is $G\cong \overline{K_n}$ or $G \cong \frac{n}{2}K_2$ when $n$ is even.
\end{theorem}

\noindent\textit{Proof.} 
From Theorem \ref{ubexvertex}, we have $\displaystyle \varepsilon_{ex_i}\leq \frac{1}{2}\left(\frac{\delta_M}{\delta_m}+\frac{\delta_m}{\delta_M}\right)\sqrt{d_i}$. Taking sum over all $i$, we have
\begin{eqnarray*}
	\varepsilon_{ex}(G)&\leq& \frac{1}{2}\left(\frac{\delta_M}{\delta_m}+\frac{\delta_m}{\delta_M}\right)\sum_{i=1}^{n}\sqrt{d_i} \nonumber\\
	&\leq& \frac{1}{2}\left(\frac{\delta_M}{\delta_m}+\frac{\delta_m}{\delta_M}\right) \left(\sum_{i=2}^{\hat{n}}\sqrt{d_i}+\sqrt{\delta_M}+\sqrt{\delta_m}\right) \nonumber\\
	&\leq& \frac{1}{2}\left(\frac{\delta_M}{\delta_m}+\frac{\delta_m}{\delta_M}\right)\left(\sqrt{(n-2)(2e-\delta_m-\delta_M)}+\sqrt{\delta_m}+\sqrt{\delta_M}\right). ~~\mbox{\Big(by Lemma \ref{poinq}\Big)}
\end{eqnarray*}

Equality holds when all the vertex degrees are equal and the connected components are isomorphic to a star graph. So, the sufficient and necessary condition for equality is either $G$ is an edgeless graph of order $n$ or $G \cong \frac{n}{2}K_2$ when $n$ is even.
\qed

\begin{remark}
    If we consider the graph $G\cong K_{1,n-1}$, then the upper bound given in Theorem \ref{ubtexenergy} is better than the upper bound in Lemma \ref{ubdexadiga} given by Adiga et al. in \cite{ar18}. So, we conclude that the bound given in Theorem \ref{ubtexenergy} is either better or incomparable with the bound given in Lemma \ref{ubdexadiga}.
\end{remark}

\begin{remark} Since $\sqrt{(n-2)(2e-\delta_m-\delta_M)}+\sqrt{\delta_m}+\sqrt{\delta_M}$ is strictly smaller than $\sqrt{2ne}$ when $\delta_m\neq\delta_M$ and both coincide when $\delta_m=\delta_M$, the bound obtained in Theorem \ref{ubtexenergy} above is an improvement of that in Das et al. (Lemma \ref{dg17ubd} (i)).
\end{remark}

\begin{remark}
The upper bound given in Theorem \ref{ubtexenergy} is incomparable with the bound given by Liu et al. in Lemma \ref{ubexlm05}. If we consider the graph $G\cong K_{1,4}$, by applying Theorem \ref{ubtexenergy}, the upper bound of extended energy is $12.75$, whereas the upper bound obtained by applying Lemma \ref{ubexlm05} is $13.046$. Evidently, the upper bound given in Theorem \ref{ubtexenergy} is better than the bound given in Lemma \ref{ubexlm05}.
\end{remark}

\begin{theorem}\label{ubdexfth}
	Let $G$ be an $(n,e)$-graph with smallest degree $\delta_m$. Also, let $F$ be the forgotten index of $G$. Then $$\varepsilon_{ex}(G)\leq \sqrt{\frac{nF}{2\delta_m^2}+ne}.$$ 
\end{theorem}

\noindent\textit{Proof.} 
From Theorem \ref{ubdvertexfth}, we can write $\displaystyle\varepsilon_{ex_i} \leq \sqrt{ \frac{1}{4\delta_m^2}\sum_{i\sim j}d_j^2+\frac{d_i^3}{4\delta_m^2}+\frac{1}{2}d_i}.$ Taking sum over all $i$, we get

\begin{eqnarray*}
	\varepsilon_{ex}(G)&\leq& \sum_{i=1}^{n}\sqrt{\frac{1}{4\delta_m^2}\sum_{i\sim j}d_j^2+\frac{d_i^3}{4\delta_m^2}+\frac{1}{2}d_i}\nonumber\\
	&\leq&\sqrt{n\sum_{i=1}^{n}\left(\frac{1}{4\delta_m^2}\sum_{i\sim j}d_j^2+\frac{d_i^3}{4\delta_m^2}+\frac{1}{2}d_i\right)} \nonumber ~~\mbox{\Big(from Lemma \ref{poinq}\Big)}\\ 
	&=&\sqrt{n\left(\frac{1}{4\delta_m^2}\sum_{i=1}^{n}\sum_{i\sim j}d_j^2+\sum_{i=1}^{n}\frac{d_i^3}{4\delta_m^2}+\frac{1}{2}\sum_{i=1}^{n}d_i\right)}\nonumber\\
	&=&\sqrt{\frac{nF}{4\delta_m^2}+\frac{nF}{4\delta_m^2}+ne}\nonumber\\
	&=&\sqrt{\frac{nF}{2\delta_m^2}+ne}.
\end{eqnarray*}
\qed

\begin{remark}
Since $\frac{\delta_M}{\delta_m}+\frac{\delta_m}{\delta_M}\geq 2$ and $F\geq 2e\delta_m^2$, the above upper bound of extended energy is better than that of Das et al. (Lemma \ref{dg17ubd} (ii)).
\end{remark}

\section{Lower bound of extended energy}
In this section, we provide some lower bounds of extended energy.

The following lower bound of extended energy is given by Das et al. \cite{dg17}.
\begin{lemma}\cite{dg17}\label{lbdexdg}
    Let $G$ be a $(n,e)$-graph. Then $$\varepsilon_{ex}(G)\geq 2\sqrt{e}$$ with equality if and only if $G\cong nK_1 $ or $G\cong K_{r,r}\cup(n-2r)K_1$, 
    $r=1,2, \dots, \lfloor{\frac{n}{2}}\rfloor$.
\end{lemma}

\begin{theorem}\label{lbdexth02}
    Let $G$ be a $(n,e)$-graph, where $e>0$. If $\varepsilon_{ex}(G)$ is the extended energy of $G$ then
    $$\varepsilon_{ex}(G)\geq \frac{2e}{k}$$ 
    where $k=\frac{1}{2}\left(\frac{\delta_M}{\delta_m}+\frac{\delta_m}{\delta_M}\right)\delta_M$ and the sufficient and necessary condition for equality being $G \cong K_{\delta_M,\delta_M}$.
\end{theorem}

\noindent\textit{Proof.}
   From Theorem \ref{lbdvexe01}, it follows.
\qed

\begin{remark}
It may be noted that the lower bound obtained through Theorem \ref{lbdexth02} sometimes yields a better lower bound than that obtained through Lemma \ref{lbdexdg}. It can be substantiated by considering the cycle graph $C_n$ for $n>4$.
\end{remark}

    \begin{theorem}\label{bldexen03}
        Let $G$ be a $(n,e)$-graph with $\det(A_{ex})\neq 0$. If $\varepsilon_{ex}(G)$ is the extended energy and $EA(G)$ is the extended adjacency index, then 
        $$\varepsilon_{ex}(G)\geq \frac{2EA(G)}{n}+n-1+\ln|\det(A_{ex})|-\ln\frac{2EA(G)}{n}$$
        with the necessary and sufficient condition for equality is $G \cong K_n$.
    \end{theorem}
    
\noindent\textit{Proof.}
Since $\det(A_{ex})\neq 0$, then $|\eta_i|>0$ for $i=1, \dots, n$. It is trivial to show that for all $x>0$ the following relation holds.
\begin{eqnarray}\label{releq01}
    x\geq 1+\ln x
\end{eqnarray}
with equality holds only when $x=1$. 

Now,
\begin{eqnarray}
    \varepsilon_{ex}(G)=\sum_{i=1}^{n}|\eta_i|&=&\eta_1+\sum_{i=2}^{n}|\eta_i|\nonumber\\
    &\geq& \eta_1 + n-1 + \sum_{i=2}^{n}\ln|\eta_i| ~~ \text{ [from equation (\ref{releq01})]}\label{equality}\\
    &=& \eta_1 + n-1 + \ln|\det(A_{ex})|-\ln\eta_1.\nonumber
\end{eqnarray}

If we consider a function $p(y)=y+n-1+ \ln|\det(A_{ex})|-\ln y$, then $p(y)$ is increasing in the interval $y\geq 1$. Since $\eta_1\geq \frac{2}{n}\sum_{i\sim j}\frac{1}{2}\left(\frac{d_i}{d_j}+\frac{d_j}{d_i}\right)=\frac{2EA(G)}{n}\geq 1$, so $p(\eta_1)\geq p\left( \frac{2EA(G)}{n}\right)$, which implies
\begin{eqnarray}
    \varepsilon_{ex}(G) &\geq& \eta_1 + n-1 + \ln|\det(A_{ex})|-\ln\eta_1=p(\eta_1)\geq p\left( \frac{2EA(G)}{n}\right)\nonumber\\
    &\geq& \frac{2EA(G)}{n}+n-1+\ln|\det(A_{ex})|-\ln\frac{2EA(G)}{n}\label{equa2} \cdot
\end{eqnarray}

Clearly, for $G\cong K_n$, equality holds. Again, the relation (\ref{equality}) will be an equality only when $|\eta_2|=\cdots=|\eta_n|=1$ and the relation (\ref{equa2}) will hold with equality only when $\eta_1=\frac{2EA(G)}{n}\cdot$ Since the trace of $A_{ex}(G)$ is zero, equality is possible only if $\eta_2=\cdots=\eta_n=-1$ and $\eta_1=n-1=\frac{2EA(G)}{n}\cdot$ Thus we conclude that the equality occurs only when $G\cong K_n$.
\qed

\begin{remark}
    Since the lower bound given in Theorem \ref{bldexen03} above is attained for all $K_n$, whereas the lower bound given by Das et al. \cite{dg17} in Lemma \ref{lbdexdg} is attained for $K_n$ only when $k=2$. Thus, for $K_n$ with $n\geq 3$, the lower bound given in Theorem \ref{bldexen03} is better than that obtained from Lemma \ref{lbdexdg}. That is the lower bound given in Theorem \ref{bldexen03} is at least incomparable with that obtained from Lemma \ref{lbdexdg}, if not better always.
\end{remark}

\section{Nordhaus-Gaddum-type bounds of extended energy}
In this section, we establish some upper as well as lower Nordhaus-Gaddum-type (N-G-type) bounds of a graph $G$, and we also show that our bounds are stronger than the existing ones. Wang et al. established the following (N-G-type) result for extended energy in \cite{wm19}.
\begin{lemma}\cite{wm19}\label{ubexnorgod1}
	Let us consider $G$ be a connected  graph of order $n$, size $e$, smallest degree $\delta_m$, largest degree $\delta_M$ and $\bar{G}$ be the complement $G$. If there are no isolated vertices in $G$ and $\bar{G}$  and $\varepsilon_{ex}(G)\geq \varepsilon_{ex}(\bar{G})$, then $$\varepsilon_{ex}(G)+\varepsilon_{ex}(\bar{G})\leq\sqrt{2ne}\sqrt{\frac{\delta_M^2(\hat{n}-\delta_m)^2}{\delta_m^2(\hat{n}-\delta_M)^2}+\frac{\delta_m^2(\hat{n}-\delta_M)^2}{\delta_M^2(\hat{n}-\delta_m)^2}+2}.$$
\end{lemma}

In the following, we obtain a new result of a similar type.
\begin{theorem}\label{ubexnorgod}
	Consider a graph $G$ with order $n$, size $e$, largest degree $\delta_M$, smallest degree $\delta_m$, and $\bar{G}$ be its complement. If there is no isolated vertices in $G$ and $\bar{G}$, then $$\varepsilon_{ex}(G)+\varepsilon_{ex}(\bar{G})\leq \frac{1}{2} \left(\frac{\delta_M}{\delta_m}+\frac{\delta_m}{\delta_M}\right)\sqrt{2ne}+\frac{1}{2}\left(\frac{\hat{n}-\delta_m}{\hat{n}-\delta_M}+\frac{\hat{n}-\delta_M}{\hat{n}-\delta_m}\right)\sqrt{n^2\hat{n}-2ne}.$$
\end{theorem}

\noindent\textit{Proof.} 
	Theorem \ref{mcbd} gives
	\begin{eqnarray}\label{ubexennor1}
		\varepsilon_{ex}(G) \leq \frac{1}{2}\left(\frac{\delta_M}{\delta_m}+\frac{\delta_m}{\delta_M}\right)\sqrt{2ne}.
	\end{eqnarray}
	
The maximum and smallest degree in $\bar{G}$ are respectively $\hat{n}-\delta_m$ and $\hat{n}-\delta_M$. Since $\bar{G}$ has no zero degree vertices,  $\hat{n}-\delta_M$ and  $\hat{n}-\delta_m$ are nonzero. Thus applying Theorem \ref{mcbd} for $\bar{G}$, one can have
	
	\begin{eqnarray}\label{ubexennor2}
		\varepsilon_{ex}(\bar{G})\leq\frac{1}{2}\left(\frac{\hat{n}-\delta_m}{\hat{n}-\delta_M}+\frac{\hat{n}-\delta_M}{\hat{n}-\delta_m}\right)\sqrt{(n\hat{n}-2e)n}.
	\end{eqnarray}

Combining inequalities (\ref{ubexennor1}) and (\ref{ubexennor2}), we have $$\varepsilon_{ex}(G)+\varepsilon_{ex}(\bar{G})\leq \frac{1}{2} \left(\frac{\delta_M}{\delta_m}+\frac{\delta_m}{\delta_M}\right)\sqrt{2ne}+\frac{1}{2}\left(\frac{\hat{n}-\delta_m}{\hat{n}-\delta_M}+\frac{\hat{n}-\delta_M}{\hat{n}-\delta_m}\right)\sqrt{n^2\hat{n}-2ne}.$$
\qed

The following theorem holds if $\varepsilon_{ex}(G)\geq \varepsilon_{ex}(\bar{G})$.

\begin{theorem}\label{ubexnorgod2}
	Let $G$ be a graph of order $n$, size $e$, largest degree $\delta_M$, smallest degree $\delta_m$  and $\bar{G}$ be the complement of the graph $G$. If $\varepsilon_{ex}(G)\geq \varepsilon_{ex}(\bar{G})$, then $$\varepsilon_{ex}(G)+\varepsilon_{ex}(\bar{G})\leq\left(\frac{\delta_M}{\delta_m}+\frac{\delta_m}{\delta_M}\right)\sqrt{2ne}$$  provided that both the graph $\bar{G}$ and $G$ are connected.
\end{theorem}

\noindent\textit{Proof.} 
If $\varepsilon_{ex}(G)\geq \varepsilon_{ex}(\bar{G})$ then we can write, 
	\begin{eqnarray*}
		 \varepsilon_{ex}(G)+ \varepsilon_{ex}(\bar{G})&\leq& 2 \varepsilon_{ex}(G)\\
		 &\le&\left(\frac{\delta_M}{\delta_m}+\frac{\delta_m}{\delta_M}\right)\sqrt{2ne}.~[\mbox{by Theorem \ref{mcbd}}]
	\end{eqnarray*}
\qed

\begin{remark}
	The N-G-type upper bound of extended energy obtained in Theorem \ref{ubexnorgod2} is stronger than that of Wang et al. (Lemma \ref{ubexnorgod1}). To show this, we claim that
	 $$\frac{\delta_M}{\delta_m}+\frac{\delta_m}{\delta_M}\leq\frac{\delta_M(\hat{n}-\delta_m)}{\delta_m(\hat{n}-\delta_M)}+\frac{\delta_m(\hat{n}-\delta_M)}{\delta_M(\hat{n}-\delta_m)} =\sqrt{\frac{\delta_M^2(\hat{n}-\delta_m)^2}{\delta_m^2(\hat{n}-\delta_M)^2}+\frac{\delta_m^2(\hat{n}-\delta_M)^2}{\delta_M^2(\hat{n}-\delta_m)^2}+2}.$$  
	 Because, otherwise
	\begin{eqnarray*}
	\frac{\delta_M(\hat{n}-\delta_m)}{\delta_m(\hat{n}-\delta_M)}+\frac{\delta_m(\hat{n}-\delta_M)}{\delta_M(\hat{n}-\delta_m)}&<&\frac{\delta_M}{\delta_m}+\frac{\delta_m}{\delta_M}\\
	\Rightarrow \frac{\delta_M}{\delta_m}\left(\frac{\hat{n}-\delta_m}{\hat{n}-\delta_M}-1\right)+\frac{\delta_m}{\delta_M}\left(\frac{\hat{n}-\delta_M}{\hat{n}-\delta_m}-1\right)&<&0\\
	\Rightarrow\frac{(\delta_M-\delta_m)(\delta_M^2-\delta_m^2)\hat{n}}{\delta_M\delta_m(\hat{n}-\delta_M)(\hat{n}-\delta_m)}&<& 0,
	\end{eqnarray*}
which is impossible.
\end{remark}

In \cite{wm19}, Wang et al. have also obtained the following N-G-type result on spectral radius of the extended adjacency matrix.

\begin{lemma}\cite{wm19}\label{lbdngsr}
	Let $G$ be an $(n,e)$-graph with largest (resp. smallest) degree $\delta_M$ (resp. $\delta_m$). If $\eta_1$ and $\bar{\eta_1}$ denote the spectral radii of $A_{ex}(G)$ and $A_{ex}(\bar{G})$ respectively, then $$\eta_1+\bar{\eta_1}\geq \frac{F}{n\delta_M^2}+\frac{\bar{F}}{n(\hat{n}-\delta_m)^2}$$ where equality holds for regular graphs.
\end{lemma}

We provide an improvement of the above bound in the following theorem.

\begin{theorem}\label{lbdngsr1}
	Let $G$ be a graph of order $n$ and $\bar{G}$ be its complement. Further, suppose that $G$ and $\bar{G}$ both are connected. If $\eta_1$ and $\bar{\eta_1}$ are the spectral radii of $A_{ex}(G)$ and $A_{ex}(\bar{G})$ respectively, then $$\eta_1+\bar{\eta_1}\geq \hat{n}.$$ 
\end{theorem}

\noindent\textit{Proof.} Let $\lambda_1$ and $\bar{\lambda_1}$ be the spectral radii of $A(G)$ and $A(\bar{G})$ respectively. Then from Lemma \ref{lbdsprth}, we have $\eta_1\geq\lambda_1\geq \frac{2e}{n}$ and $\bar{\eta_1}\geq \bar{\lambda_1}\geq \frac{n\hat{n}-2e}{n}\cdot$ Utilizing both of these, we have $$\eta_1+\bar{\eta_1}\geq \hat{n}.$$ \qed	

\begin{remark}
	The lower bound of $\eta_1+\bar{\eta_1}$ which is given in Theorem \ref{lbdngsr1} is better than that given by Wang et al. (Lemma \ref{lbdngsr}). Because otherwise, 
	\begin{eqnarray*}
		\hat{n}&<&\frac{F}{n\delta_M^2}+\frac{\bar{F}}{n(\hat{n}-\delta_m)^2}\\
		\Rightarrow n\hat{n}&<&\sum_{r=1}^{n}\frac{d_r^3}{\delta_M^2}+\sum_{r=1}^{n}\frac{\bar{d_r}^3}{(\hat{n}-\delta_m)^2}\\
		&=&\sum_{r=1}^{n}d_r\frac{d_r^2}{\delta_M^2}+\sum_{r=1}^{n}(\hat{n}-d_r)\frac{(\hat{n}-d_r)^2}{(\hat{n}-\delta_m)^2}\\
		&\leq&\sum_{r=1}^{n}d_r+\sum_{r=1}^{n}(\hat{n}-d_r)\\
		&=&n\hat{n}
	\end{eqnarray*}
which is impossible.
\end{remark}

Finally, in this section, we find an N-G-type lower bound of the extended energy of some connected graphs.
\begin{theorem}\label{norgad3}
	Let $G$ and its complement $\bar{G}$ be both connected graphs of order $n$. Then the following relation holds
    $$\varepsilon_{ex}(G)+ \varepsilon_{ex}(\bar{G})\geq 2\hat{n}.$$
\end{theorem}

\noindent\textit{Proof.} 
	Let the spectral radii of $A_{ex}(G)$ and $A_{ex}(\bar{G})$ are $\eta_1$ and $\bar{\eta_1}$ respectively. Then from Theorem \ref{lbdngsr1}, we can write
	\begin{eqnarray*}
		\varepsilon_{ex}(G)+\varepsilon_{ex}(\bar{G})&\geq& 2\eta_1+2\bar{\eta_1}\\
		&\ge& 2\hat{n}.
	\end{eqnarray*}
\qed

\begin{remark}
The bound given in Theorem \ref{norgad3} was shown to be true for regular graphs by Wang et al. \cite{wm19}. But, Theorem \ref{norgad3} covers the non-regular graphs as well.
\end{remark}

\section*{Acknowledgement}
The first author sincerely acknowledges the financial assistance of UGC, India in the form of a senior research scholarship to carry out this work.

\end{document}